\documentclass[11pt]{article}%
\usepackage{amsmath,amssymb,enumerate}%
\usepackage{amsmath}%
\usepackage{amsfonts}%
\usepackage{amssymb}%
\usepackage{graphicx}
\newtheorem{definition}{Definition}
\newcommand{\tmop}[1]{\ensuremath{\operatorname{#1}}}
\newtheorem{theorem}{Theorem}
\newcommand{\tmem}[1]{{\em #1\/}}

\newenvironment{proof}{
\noindent\textbf{Proof}\ }{\hspace*{\fill}
\begin{math}\Box\end{math}\medskip}
\newenvironment{enumeratenumeric}{\begin{enumerate}[1.]}{\end{enumerate}}
\newtheorem{lemma}{Lemma}
\newtheorem{varremark}{Remark}

\newtheorem{proposition}{Proposition}
\newcommand{\tmstrong}[1]{\textbf{#1}}

\begin{document}

\title{Random walks on the hypergroup\\of circles in a finite field}
\author{Le Anh Vinh\\School of Mathematics\\University of New South Wales\\Sydney
2052 Australia}
\maketitle

\begin{abstract}
In this paper we study random walks on the hypergroup of circles in a finite
field of prime order $p=4l+3$. We investigate the behavior of random walks on
this hypergroup, the equilibrium distribution and the mixing times. We use two
different approaches---comparison of Dirichlet Forms (geometric bound of
eigenvalues), and coupling methods, to show that the mixing time of random
walks on hypergroup of circles is only linear.

\end{abstract}

\thispagestyle{empty} \thispagestyle{empty}

\section{Preliminary}

Suppose that $p$ is a prime of the form $p=4l+3$ for some integer $l,$ and
that $\mathbb{F}_{p}$ is the field with $p$ elements. The following definition
follows [7], where the importance of the notion of
quadrance is developed.

\begin{definition}
The \textbf{circle} in $F_{p}$ with \textbf{center} $A_{0}=[x,y]$ and
\textbf{quadrance} $K$ is the set of all points $\left[  u,v\right]  $ in
$F_{p}\times F_{p}$ such that%
\[
(u-x)^{2}+(v-y)^{2}=K\left(  \operatorname{mod}p\right)  .
\]

\end{definition}

There are exactly $p$ circles with center $O=[0,0],$ with quadrances
$0,1,2,\ldots,p-1$; we denote them by $C_{0},C_{1},C_{2},\ldots,C_{p-1}$
respectively. If we start from $O=[0,0]$, take a random step by translating by
an element of $C_{i}$, and then take another random step by translating by an
element of $C_{j}$, the final point will be an element of $C_{k}$ for some
$k.$ Counting over all possible such combinations, there are $N_{ij}^{k} $
ways to reach to a point of $C_{k}$ by using steps from $C_{i}%
,\operatorname{then}C_{j}$ randomly. We can write this relation as
\[
C_{i}C_{j}=\sum_{k}N_{ij}^{k}C_{k},
\]
where $N_{ij}^{k}$are non-negative integers. Let $n_{ij}^{k}=\frac{N_{ij}^{k}%
}{|C_{i}\Vert C_{j}|}$ then it can be written as distribution form%

\[
c_{i} c_{j} = \sum_{k} n_{i j}^{k} c_{k}%
\]

where $n_{i j}^{k} \geqslant0$ and $\sum_{k} n_{i j}^{k} = 1$ for any $i, j$.

\begin{theorem}
Suppose $p$ is a prime of the form $p=4l+3$ for some integer $l$. Given $i,j$
non-zero and $k$ in $\mathbb{F}_{p},$ let $V=ij-\frac{(k-i-j)^{2}}{4}.$ Then
we have the following formulas for $N_{ij}^{k}$ and $n_{ij}^{k}$. i) If $V$ is
not a square in $F_{p}$ then $N_{ij}^{k}=0$ and $n_{ij}^{k}=0$. ii) If $V=0$
then $N_{ij}^{k}=|C_{i}|=p+1$ and $n_{ij}^{k}=\frac{1}{p+1}$. iii) Otherwise,
$N_{ij}^{k}=2|C_{i}|=2(p+1)$ and $n_{ij}^{k}=\frac{2}{p+1}$.
\end{theorem}

\begin{proof}
  In this case $|C_{0}|=1,|C_{i}|=p+1\operatorname{for}i=1,\ldots,p-1$. We start
from $[0,0]$, and translate by $C_{i}$, to reach to a point, say $(x,y)$ where
$x^{2}+y^{2}=i$ ($\operatorname{mod}p$) (we will not write mod $p$ from now).
From $(x,y)$, we translate by $(u,v)$ in $C_{j}$ so that $u^{2}+v^{2}=j$
(there are $|C_{j}|$ such possible moves from $(x,y)$). We will reach a point
in $C_{k}$ if and only if
\[
(x+u)^{2}+(y+v)^{2}=k
\]
or
\[
xu+yv=\frac{k-i-j}{2}
\]
Recall the well-known identity
\[
(x^{2}+y^{2})(u^{2}+v^{2})=(xu+yv)^{2}+(yu-xv)^{2}.
\]
Thus the condition is equivalent to $(yu-xv)^{2}=ij-\frac{(k-i-j)^{2}}{4}=V$.

  (i) If $V$ is not a square in $F_p$ then it is clear that there does not exist such
  $x, y, u, v$. Hence $N_{i j}^k = n_{i j}^k = 0$. 
  
  (ii) If $V = 0$ or $i j = ( k - i - j )^2/4$ then $( y u - x v )^2 = 0$ or $y u = x v$. We can assume that $x \neq 0$. From $v = y u/x$, we have $j = u^2 + v^2 = ( ( y/x )^2 + 1 ) u^2 = i u^2/ x^2$ or $u^2/x^2 = j/i$. Since $i j$ is square
  so $j/i = m^2$ for some $1 \leqslant m \leqslant \frac{p
  - 1}{2}$. Thus $u/x = \pm m$, $u = \pm m x$, and $v = \pm m
  y$.
  
  We need $x u + y v = (k - i - j)/2$, so $\pm m ( x^2 + y^2 ) =
  (k - i - j)/2$ or $\pm \sqrt{i j} = (k - i - j)/2$. Since
  $( k - i - j)^2/4 = i j$ we can choose exactly either $m$ or $- m$
  to make this requirement hold. Therefore, for each $( x, y )$ in
  $C_i$, we have only one $( u, v )$ in $C_j$ such that if we go by $( x,
  y )$ followed by $( u, v )$, we reach to a point in circle $C_k$. In other
  words, $N_{i j}^k = |C_i | = p + 1$.
  
  (iii) If $V = m^2 \neq 0$. Set $a = (k - i - j)/2$ and $b = y u - x v$ ( $b = m$ or $b =    - m$), it follows that $x u + y v = a$ and $y u - x v = b$. Solving these equations we      have $u = (a x + b y)/i$ and $v = (a y - b x)/i$. It is straight
  forward to verify that $( u, v ) \in C_j$ and $( x_{} + u, y + v ) \in
  C_k$ if $( x, y ) \in C_i$. Since $m \neq 0$, for each $( x, y ) \in
  C_i$, we have two posibilities to go to $C_k $ by $C_j$. Therefore $N_{i j}^k
  = 2| C_i | = 2 ( p + 1 )$, completing the proof.
\end{proof}

To see the fact that the circles equipped with random walks product create a
hypergroup structure, recall the formal definition of (general)
hypergroup (see [6]).

\begin{definition}
  A (finite) general hypergroup is a pair $(\mathcal{K},\mathcal{A})$ where
  $\mathcal{A}$ is a *-algebra with unit $c_0$ over $C$ and
  $\mathcal{K}= \{ c_0, c_1, \ldots, c_n \}$ is a subset of $\mathcal{A}$
  satisfying
  \begin{enumeratenumeric}
    \item $\mathcal{K}$ is a bsis of $\mathcal{A}$    
    \item $\mathcal{K}^{\ast} =\mathcal{K}$    
    \item The struture constants $n_{i j}^k \in C$ defined by    
     \[ c_i c_j = \sum_k n_{i j}^k c_k \]
    satisfy the conditions    
     \[c_i^{\ast} = c_j \Leftrightarrow n_{i j}^0 > 0, \]
     \[c_i^{\ast} = c_j \Leftrightarrow n_{i j}^0 > 0 \]      
  \end{enumeratenumeric}
\end{definition}

$\mathcal{K}$ is called {\tmem{hermitian}} if $c_i^{\ast} = c_i$ for all $i$,
commutative if $c_i c_j = c_j c_i$ for all $i, j$, {\tmem{real}} if $n_{i
j}^k$ if $n_{i j}^k$ $\in R$ for all $i, j, k$, {\tmem{positive }}if
$n_{i j}^k \geqslant 0$ for all $i, j, k$ and {\tmem{normalized }}if $\sum_k
n_{i j}^k = 1$ for all $i, j$. A generalized hypergroup which is both positive
and normalized will be called a hypergroup. There are board examples and
applications of (generalized) hypergroups which can be found in [7].

Now, we show that the set $C = \{ c_0, c_1, \ldots, c_{p - 1} \}$ with the
relation $c_i c_j = \sum_{k = 0}^{p - 1} n_{i j}^k c_k$ (defined above) is
indeed a hypergroup. Let $n_{i j}^k = \frac{N_{i j}^k}{|C_i \| C_j |}$
then $c_i c_j = \sum_k n_{i j}^k c_k $. It is clear that $n_{i j}^k \geqslant
0$, and $\sum_k n_{i j}^k = 1$ for any $i, j$. From Theorem
1, $n_{i j}^0 \neq 0$ if and only if $i j - ( 0 - i - j )^2/4 = - ( i - j )^2/4$ is square. If $p \equiv 3$ (mod 4) then $-1$ is not a square in $F_p$, so $n_{i j}^0 \neq 0$ if and only if $i = j$. Let $c_i^{\ast} = c_i$ then $C$ is a hermitian commutative hypergroup
(note that, $n_{i j}^k$ is symmetric w.r.t $i, j$ and $k$ so $C$ is commutative).
For $p \equiv 1$ (mod $4$), the null circle is extraodinary and yields
many troubles in calculations so we restrict our interests in case $p \equiv 3$
(mod $4$) throughout this paper.

We also need the following theorem which is first proved by the author in [6] (a long proof using quadratic residue method). We present here a robust proof obtained directly
from Theorem 1.

\begin{theorem}
  For any $i, j, k \in F_p$ ($p \equiv 3$ mod $4$) , we write $c_i c_j c_k =
  \sum_l \alpha_{i j k}^l c_l$ then $\alpha_{i j k}^l > 0$ if $i, j, k, l$ are
  all non-zero. If one of $i, j, k, l$ is zero, let $x, y, z$ be three
  remaining numbers then $\alpha_{i j k}^l > 0$ iff $x y - ( z - x - y)^2/4$ is square in $F_p$.
\end{theorem}

\begin{proof}
  If one of $i, j, k, l$ is zero then it becomes: for three remaining numbers, $x, y$,and     $z$ are we have a triangle with quadrances $x, y$ and $z$ (note
  that a step from $C_0$ when $p \equiv 3$ (mod $4$) is a stationary
  step). By Theorem 1, it holds if and only if $V = x y -
  ( z - x - y )^2/4$ is square. Now we assume that $i,
  j, k$ and $l$ are all nonzero. It is clear that $C_l$ is reachable by steps from
  $C_i, C_j$ then $C_k$ if and only if we can start from $( 0, 0 )$, go by
  steps from $C_i, C_j, C_k$ then $C_l$ to come back to $( 0, 0 )$. Now,
  from $( 0, 0 )$, go a step from $C_i$, followed by a step from $C_j$ we have $|C_i \| C_j|
  = ( p + 1 )^2$ possible steps. From Theorem 1, there is no more than $2 ( p
  + 1 )$ steps that can reach the same circles. So, by choosing steps from
  $C_i$ then $C_j$, we can reach at least $( p + 1 )^2/2 ( p + 1 )
  = (p + 1)/2$ circles.  Applying the same argument, by choosing steps from
  $C_l$ then $C_k$ (start from $( 0, 0 )$), we can reach at least $(p + 1)/2$ circles. Since we have only $p$ circles, by the Pigeon Holes
  Principle, there exists a circle say $C_t$ that is reachable from both directions.
  Therefore, we can go by $C_i, C_j$ from $( 0,0 )$ to $C_t$; then from $C_t$,
  go by $C_k$ to $C_l$. The statement follows.
\end{proof}

\section{Random walks on hypergroup of circles}

In this section, we will consider the random walk by $C_1$, i.e. we choose all steps from the unit circle $C_1$. We call it random
walk $C_1$ by abuse of notation. In general, at $n^{\tmop{th}}$ step we have the
relation
\[ c_1^n = \sum_{_{j = 0}^{}}^{p - 1} \alpha_{n, j} c_j \]
where $\alpha_{n, j} \geqslant 0$ for $j = 0, \ldots, p - 1$ and $\sum_{j = 0}^{p - 1} \alpha_{n,j}$ = 1.

From Theorem 2, after no more than $4$ steps, we will reach all
circles at every state. Thus, $\alpha_{n, j} > 0$ for all $i = 0, \ldots, p - 1$
and $n \geqslant 4$.

Let $K$ be a Markov kernel. The probability $\pi$ is invariant or stationary
for $K$ if $\pi K = \pi$. A Markov kernel $K$ is irreducible if for any two
states $x, y$ there exists an integer $n = n ( x, y )$ such that $K_n ( x, y )
> 0$. A state $x$ is called aperiodic if $K_n ( x, x ) > 0$ for all
sufficiently large $n$. If $K$ is irreducible and has an aperiodic state then
all states are aperiodic and $K$ is {\tmem{erogodic}}.

Let $K$ be $c_1$, i.e $K ( c_i, c_j ) = n_{i 1}^j$. Our main interest in
this section is the iterated kernel $K_n ( x, y )$. We write $K$ as
$c_1$ and $K_n$ as $c_1^n$. By Theorem 2, from $C_j$ we can go to $C_i$
by no more than 2 steps from $C_1$. Thus, $K$ is irreducible. Besides,
$c_1^n ( c_1, c_1 ) > 0$ for all $n > 3$, so $K$ has an aperiodic state. Thus
all states are aperiodic.

The following definition gives us the total variation distance between two
probability measures.

\begin{definition}
  Let $\mu, \nu$ be two probability measures on the set $X$. The total
  variation distance is defined by  
  \[ d_{\tmop{TV}} ( \mu, \nu ) = \| \mu - \nu \|_{\tmop{TV}} = \max_{_{A
  \subset X}} | \mu ( A ) - \nu ( A ) | = \frac{1}{2} \sum_{x \in X} | \mu ( x
  ) - \nu ( x ) |. \]
\end{definition}

Ergodic Markov chains are useful algorithmic tools in which, regardless of
their initial state, they eventually reach a unique stationary distribution.
The following theorem, originally proved by Doeblin, details the essential
property of ergodic Markov chains.

\begin{theorem}
  Let $K$ be any ergodic Markov kernel on a finite state space $X$ then $K$
  admits a unique stationary distribution $\pi$ such that
 \[
  \forall x, y \in X, \lim_{t \rightarrow
  \infty} K_t ( x, y ) = \pi ( y ). \]

\end{theorem}

From Theorem 3, the Markov kernel $K$ (i.e $C_1$)  admits a unique invariant
distribution $\pi$. For all $c_i, c_j$ $\lim_{n \rightarrow
\infty} c_1^n ( c_i, c_j ) = \pi ( y )$. In general, it is difficult to
determine this unique stationary distribution. However, in this case, it can
be found easily by the following lemma.

\begin{lemma}
  If there exists a distribution $\pi$ such that
  
\begin{center}
    $\pi ( x ) K ( x, y ) = \pi ( y  ) K ( y, x )$ for all $x, y \in X$
\end{center}
  
  then such a $\pi$ is a stationary distribution. (in this case, $K$ is called
  {\tmem{reversible}})
\end{lemma}

From Theorem 3 and Lemma 1, we now can tell exactly the behaviour of $c_1^n$
as $n$ goes to infinity.                  

\begin{theorem}
  Over the finite field $F_p ( p \in \mathcal{P}> 3 )$, then
  \[ \lim_{n \rightarrow \infty} c_1^n = \frac{1}{p^2} c_0 + \frac{p + 1}{p^2} ( c_1 + c_2 + \ldots + c_{p - 1} ).\]  
\end{theorem}

\begin{proof}
  We know that the kernel Markov $A_1$ (of $c_1$) is erogodic. From Theorem 3, there
  exists a unique stationary distribution $\pi$ of $A_1$. Let's consider a
  distribution $\sigma$ on $C$ which is defined as $\sigma ( c_0 ) = 1/p^2$ and
  $\sigma ( c_j ) = (p + 1)/p^2$ for all $j > 0$. We see that $A_1 (
  c_i, c_j ) = c_1 ( c_i, c_j ) = n_{i 1}^j$. By Lemma 1, to show that $\pi = \sigma$, we only need to verify that $\sigma ( c_i ) n_{i 1}^j = \sigma ( c_j )
  n_{j 1}^i$. If $i, j \neq 0$ then this equation clearly holds by the
  symmetricity of $n_{i j}^k$ w.r.t $i, j$ and $k$ (when $i, j, k \neq 0$). If
  $i = 0$, then $n_{01}^j, n_{j 1}^0 \neq 0$ if and only if $j = 1$. And if $j
  = 1$ then $n_{01}^1 = 1$ and $n_{11}^0 = 1/(p + 1)$. The equation
  still holds. Thus, for any $i$ and $j$, we have $\sigma ( c_i ) n_{i 1}^j =
  \sigma ( c_j ) n_{j 1}^i$. Therefore, $\pi = \sigma$, completing the proof.
\end{proof}

Note that $|C_{o}|=1,|C_{1}|=|C_{2}|=\ldots=|C_{p-1}|=p+1$, and the space has
$p^{2}$ points, so the distribution of $c_{1}^{n}$ is, in some sense, close to
uniform over the space $F_{p}^{2}$ when $n$ tends to infinite. Walking
randomly by any $C_{i}(i\neq0$) we have the same results as for $C_{1}.$ In
hypergroup language the limiting distribution is the Haar measure on the hypergroup.

We proved that $c_1^n$ tends to the unique stationary distribution of the
hypergroup of circles, but we give no infomation about the rate of
convergence as a function of the size of the hypergroup, i.e $|C| =
p$. We define the {\tmem{mixing time}} $\tau_p ( \varepsilon )$ as
the time until the chain is within variation distance $\varepsilon$ from the
worst initial state. We give a formal definition for this concept.

\begin{definition}
  $\tau_p ( \varepsilon ) = \max_{i \in \{ 0, \ldots, p - 1 \}} \min \{ t :
  d_{\tmop{TV}} ( c_1^t ( c_i, . ), \pi ) \leqslant \varepsilon \}$ where
  $\pi$ is the distribution in Theorem $4$.
\end{definition}

We can fix $\varepsilon$ as any small constant. A popular choice is to set
$\varepsilon = 1/2 e$. We then boost to arbitrary small variation
distance by the following lemma.

\begin{lemma}
  $\tau_p ( \varepsilon ) \leqslant \tau_p ( 1/2 e ) \ln ( 1/e )$.
\end{lemma}

We now want to estimate $\tau_p = \tau_p ( 1/2 e )$. We will present
two different approaches to estimate this mixing time.

\subsection{Comparision of Dirichlet Forms}

The first approach is using the comparsion of Dirichlet forms to estimate
$\tau_p$. We first need some preliminaries. Let $X$ be a finite set. Let $K (
x, y ), \pi ( x )$ be a reversible irreducible Markov chain on $X$. Let $l^2 (
X )$ have scalar product
\[ < f, g > = \sum_{x \in X} f ( x ) g ( x ) \pi ( x).\]
From the reversibility of $\pi$ and $K$, the operator $f \mapsto K f$,
with $K f ( x ) = \sum f ( y ) K ( x, y )$, is self-adjoint on $l^2$ with
eigenvalues $\alpha_1 = 1 > \alpha_2 \geqslant \ldots \geqslant \alpha_{|X| -
1} \geqslant - 1$. These eigenvalues can be characterized by the Dirichlet
form $D$, which is defined as
\[D ( f, f ) = < ( I - K ) f, f > = \frac{1}{2} \sum_{x, y} ( f
( x ) - f ( y ) )^2 \pi ( x ) K ( x, y ).\]    

Let $V$ be a subspace of $L^2 ( X )$, set
\begin{center}
$M_D ( V ) = \max \{ D ( f, f ) ; \| f \|_2 = 1, f \in V \},$
$m_D ( V ) = \min \{ D ( f, f ) ; \| f \|_2 = 1, f \in V \}.$
\end{center}

We have the minimax characterization of eigenvlaues for $0 \leqslant i
\leqslant |X| - 1$ (for details see [5]),
\[1 - \alpha_i = \min \{ M_D ( V ) | \dim V = i + 1 \} = \max \{
m_D ( V ) | \dim V^{\bot} = i \}.\] 

Suppose we have a second reversible Markov chain on $X$, say $K', \pi'$ with
eigenvalues $\alpha_{i}^{'}$. From the minimax characterization: for
$1 \leqslant i \leqslant |X| - 1$ we have: 
\begin{equation}
\alpha_{i} \leqslant 1 - \frac{a}{A} ( 1 - \alpha_{i}^{'} ) \tmop{if} 
D' \leqslant A D, \pi' \geqslant a \pi.
\end{equation}

Set $E = \{ ( x, y ) | K ( x, y ) > 0 \}$, and $E' = \{ ( x, y ) |K' ( x, y ) >
0 \}$. For each pair $x \neq y$ with  $K' ( x, y ) > 0$, we fix a sequence
of steps $x_0 = x, x_1, \ldots, x_k = y$ with $K ( x_i, x_{i + 1} ) > 0$. This
sequence of steps is called a path $\gamma_{x y}$ of length $| \gamma_{x y} |
= k$ (paths may have repeated vertices, but a given edge appear at most once in
a given path). Set $E' ( e ) = \{ ( x, y ) \in E' | e \in
\gamma_{x y} \}$ for $e \in E$. We have the following theorem for a
bound of $A$ in (1), which will give an upper bound for the second largest
eigenvalue.

\begin{theorem}
  $([2])$ Let $K, \pi$ and $K', \pi'$ be reversible Markov chains on a
  finite set $X$. For the Dirichlet forms defined as above, then $D' \leqslant A D$
  where  
  \begin{equation}
   A = \max_{( z, w ) \in E} \{ \frac{1}{\pi ( z ) K ( z, w )}
  \sum_{E' ( z, w )} | \gamma_{x y} | \pi' ( x ) K' ( x, y ) \}.
  \end{equation}
\end{theorem}

Let $\sigma_x$ be a cycle from $x$ to $x$ with an odd number of edges (again,
we may have repeated vertices but not repeated edges). For irreducible
aperiodic chains, such cycles always exist. Suppose we have a fix
collection of such cycles for all circles $x$'s (one cycle for each circle). We define the cycle length by chain as 
\[ | \sigma_x |_K = \sum_{( z, w ) \in \sigma_x}
\frac{1}{K^{\ast} ( z, w )}\]
where $K^{\ast} ( z, w ) = \pi ( z ) K ( z, w )$
(note that $K$ is reversible so $\pi ( z ) K ( z, w ) = \pi ( w ) K ( w, z
)$). We have a following lower bound for the smallest eigenvalue of an
irreducible aperiodic Markov chain.

\begin{theorem}
  $([4])$ Let $K, \pi$ be an irreducible aperiodic Markov chain, the smallest eigenvalue $\alpha_{\min} = \alpha_{|X| - 1}$ satisfies
  \begin{equation}
  \alpha_{\min} \geqslant - 1 + \frac{2}{v}
  \end{equation}  
  where $v = \max_{_e} \sum_{_{\sigma_x \ni e}} | \sigma_x |_K \pi ( x )$.
\end{theorem}

The relation between the mixing time and eigenvalues of the kernel is
explained in the following theorem.

\begin{theorem}$([4])$ $\| K_t ( x, . ) - \pi \|_{\tmop{TV}}
  \leqslant \pi_{\ast}^{- 1/2} \alpha_{\ast}^t/2$ 
  where $\pi_{\ast} = \min_x \{ \pi ( x ) \}$ and $\alpha_{\ast} = \max \{ |
  \alpha_{|X| - 1} |, \alpha_1 \}$.
\end{theorem}

We now can give an estimation for the mixing time of random walk on
hypergroup of circles in finite fields.

\begin{proposition}
  $\tau_p = O ( p \ln p )$.
\end{proposition}

\begin{proof}
  Let we explain how we obtain the result before going into the details.
  The proof contains 3 steps:
  
  {\tmstrong{I.}} Show that $\alpha_1 \leqslant 1 - c/p$ for some $c$
  by using Theorem 5.
  
  {\tmstrong{II. }}Show that $\alpha_{p - 1} \geqslant - 1 + d/p$ for
  some $d$ by using Theorem 6.
  
  {\tmstrong{III.}} Let $\eta = \min ( c, d )$ then $\alpha_{\ast}
  \leqslant 1 - \eta/p$. We then use Theorem 7 to conclude
  the proof.
  
  Now, we go into the details.
  
  {\tmstrong{I.}} Consider two reversible Markov chains, $( A_1$, $\pi )$ and
  $( A_{1}^{'}, \pi )$ where $A_1$ is the Markov kernel for our random walk, $\pi$
  is the stationary distribution in theorem 4, and $A'_1$ is the
  equilibrium Markov kernel of $A_1$. Then $A_{1}^{'}$ is the matrix with $p$
  rows of the distribution $\pi$. Therefore, all eigenvalues of $A_{1}^{'}$ are
  $1$ with multiple 1, and 0 with multiple $p - 1$. From (1), we have $a = 1,
  \alpha_{i}^{'} = 0$, so if $D' \leqslant A D$ then $\alpha_i \leqslant 1 -
  1/A$. We want to approximate
  \begin{equation}
  A = \max_{( z, w ) \in E} \{ \frac{1}{\pi ( z ) A_1 ( z, w
  )} \sum_{E' ( z, w )} | \gamma_{x y} | \pi ( x ) A_1' ( x, y ) \}.
  \end{equation}
  
  Then we have $\alpha_1 \leqslant 1 - 1/A$.
  
  First, we need to fix a path for each edge of $A_{1}^{'}$. It is clear that $(
  x, y ) \in E'$ for all $x, y$. For each $x \neq y$, $x = c_r$, $y = c_s$
  for some $r, s$. We assume that $r < s$. We define $\gamma_{x y}$ as follow:
  
  $\diamond$ If $x = c_0, y = c_1$ then $\gamma_{x y} = \{ c_o, c_1 \}$.
  
  $\diamond$ If $x = c_0, y \neq c_1$ then from Theorem 2, there exists $c_k$
  such that $c_1 ( c_1, c_k ) > 0$ and $c_1 ( c_k, y ) > 0$. We set
  $\gamma_{x y} = \{ c_0, c_1, c_k, y \}$.
  
  $\diamond$ If $x, y \neq c_0$, then from Theorem 2, there exists $c_k$ such
  that $c_1 ( x, c_k ) > 0$ and $c_1 ( c_k, y ) > 0$. We set $\gamma_{x,
  y} = \{ x, c_k, y \}$  (**)
  
  Let $( z, w )$ be any edge of $E$, then $z = c_r$, $w = c_s$ for some $r,
  s$, we assume that $r < s$. There are 2 cases of $r$.
  
  {\tmstrong{Case 1.}} $r = 0$, then $s = 1$ (by Theorem 1). All paths
  $\gamma_{x y}$ that contain this edge are of form  ${ x = c_0,c_1,c_k,y}$
  (we assume that quadrance of circle $x$ is less than quadrance of circle
  $y$). By the way we define $\gamma_{x y}$, there are $p - 1$ such paths and
  each path is of length 4. We have $\pi ( x ) = \pi ( c_0 ) = 1/p^2$ and
  $A_1' ( c_0, y ) \leqslant (p + 1)/p^2$. Thus, if $r = 0$, we have
  \[ \frac{1}{\pi ( z ) A_1 ( z, w )} \sum_{E' ( z, w )} | \gamma_{x y} | \pi (
  x ) A_1' ( x, y ) \leqslant 4 p^2 ( p - 1 ) \frac{p + 1}{p^4} \leqslant 4.\]   
  
  {\tmstrong{Case 2.}} $1 < r < s$, then all path $\gamma_{x y}$ that contain
  this edge are of form $\{ x = c_r, c_s, c_k = y \}$ or $\{ x = c_k, c_r, c_s
  = y \}$ for some $c_k$. In order to have a path like that, we need $k$ such that
  $r s - ( k - r - s )^2/4$ is nonzero square (by Theorem 1). Thus, we have no more than $( p + 3 ) / 2$ of such $k$ (note that not all $k$ are choosen in (**), and we are overcounting. However, it does not matter since
  we are considering the worst case). Therefore, there are no more than $p + 3$ of such
  paths, each path is of length 3. We have, $\pi ( x ) = (p + 1)/p^2$, $A_{1}^{'} ( x, y ) =  (p + 1)/p^2$, $\pi ( z ) = (p + 1)/p^2$ and $A_1
  ( z, w ) \geqslant 1/(p + 1)$ (from Theorem 2 and Theorem 4). Therefore,
  
  \begin{equation}
  \frac{1}{\pi ( z ) A_1 ( z, w )} \sum_{E' ( z, w )} | \gamma_{x y} | \pi (
  x ) A_1' ( x, y ) \leqslant 3 \frac{( p + 3 ) ( p + 1 )^2}{p^2}.
  \end{equation}  
  
  From (4) and (5), we have $A \leqslant 3 ( p + 3 ) ( p + 1 )^2/p^2$. Hence      \[\alpha_1 \leqslant 1 - \frac{p^2}{3 ( p + 3 ) ( p + 1 )^2} \leqslant
  1 - \frac{c}{p} \]
  for some constant $c$.                            
  
  {\tmstrong{II.}} We will construct a collection of cycles for $A_1$. We put
  $c_0$ aside and only consider $p - 1$ circles $c_1, \ldots, c_{p - 1}$.
  We call $c_i$ and $c_j$ adjacent if and only if $A_1 ( c_i, c_j ) > 0$.
  If $c_i$ is self-adjacent then we let $\sigma_{c_i} = \{ c_i \}$. From the
  remaining circles, if possible, we choose any 2 adjacent circles, say $c_s$
  and $c_t$, then $c_s$, $c_t$ are not self-adjacent. By the proof of
  Theorem 2, both of them are adjacent to at least $(p + 1)/2$ circles, so
  each of them adjacent to at least $(p - 1)/2$ circles from $C
  \backslash \{ c_s, c_t \}$ (where $C = \{ c_0, \ldots, c_{p - 1} \}$ is the
  hypergroup). By Pigeon Hole Principal, there exists a circle $c_u$
  ($u \neq s, t$ but $u$ may be zero) which is common adjacent to these
  circles. We set $\sigma_s = \sigma_t = \{ c_s, c_t, c_u \}$. If
  $\sigma_u$ has not yet determined, we also let $\sigma_u = \sigma_s$.
  
  We keep doing until we have choosen all circles or we get stuck. It is easy to
  see that by the way of choosing cycles, we have not choosen any edge twice
  (in each step, we introduce a triangle that may contain an old vertice but
  always has 3 new edges). If we have choosen all circles, we finish our
  construction. If not, let we call $U$ the set of circles that their
  cycles have been determined but not $c_0$ (even if we have defined
  $\sigma_{c_0}$ in some step), and $V$ for the rest, except $c_0$. From the
  construction, we see that all vertices in $V$ are not adjacent to one another.
  If $V$ has more than one circle, let $c_s$ and $c_t$ be any two
  circles in $V$. Then by proof of Theorem 2, $c_s$ has at least $(p -
  1)/2$ adjacent circles in $U$, and so is $c_t$. We have $|U| \leqslant p -
  3$ (since $c_0, c_s, c_t$ are not in $U$). By pigeon hole principle, there
  exists two circles $c_u, c_v$ which are common adjacent to $c_s, c_t$. Let
  $\sigma_u = \{ c_u, c_{u_1}, c_{u_2} \}$, we define $\sigma_s = \sigma_t =
  \{ c_s, c_v, c_t, c_u, c_{u_1}, c_{u_2}$,$c_u ( c_s ) \}$ (odd cycles). 
  
  Suppose we have choosen $k$ pairs in $V$, and conresponding $k$ pairs in $U$
  in such a way that they are all disjoint. Let's call the set of 2$k$ circles from $k$
  pairs in $V$ (and $U$) $V_1$ (and $U_1$). If $V \backslash
  V_1$ has more than one circle, let $c_m, c_n$ be any two circles in
  this set. By the proof of Theorem 2, $c_m$ is adjacent to at least $(p -
  1)/2 - 2 k$ circles in $U \backslash U_1$, and so is $c_n$. We have $|U
  \backslash U_1 | \leqslant p - 3 - 4 k = 2 ( (p - 1)/2 - 2 k ) + 2$
  (since $U_1, V_1$ and $c_0, c_n, c_m$ are not in $U$). Again, by Pigeon
  Hole Principle, there exists two circles $c_h, c_l \in U \backslash U_1 $
  which are common adjacent to $c_m, c_n$. Thus, we can do similarly as above for these four circles.
  
  We keep doing until we have no more than one circle left in $V$.
  During the process, each pair in $V$ and its conrespoding (in $U$) are
  disjoint so all edges from $V$ to $U$ are used in only 2 cycles. Consider
  any cycle in $U$ which is triangle by the construction. When choosing pairs
  of circles from $V$, vertices of $U$ are chosen no more than 3 times. So in
  total, any edge of any cycle in $U$ belongs to no more than $1 + 3 \times
  2 = 7$ cycles. Now, we have at most 2 circles left (at most one from $V$ and
  $c_0$). For any circle whose its cycle has not yet determined we assign to
  it any odd cycle of length $\leqslant 5$ that contains it (there exists such
  cycle by Theorem 2). So, each edge belongs to no more than 9 cycles,
  and each cycle has no more than 7 edges. 
  
  For each $( z, w )$ such that $A_1 ( z, w ) > 0$ we have 
  \[A_1^{\ast} ( z, w ) = \pi ( z ) A_1 ( z, w ) \leqslant 2/p^2.\]
  So for each
  cycle $\sigma_{c_i}$, $| \sigma_{c_i} |_{A_1} \leqslant 7 p^2$. We
  have 
  
 \[ v = \max_{_e} \sum_{_{\sigma_{c_i} \ni e}} | \sigma_{c_i} |_{A_1} \pi
  ( c_i ) \leqslant 9 \times 7 p^2 \times \frac{p + 1}{p^2} = 63 ( p + 1 ). \] 
  
  Hence, by Theorem 6 we have $\alpha_{p - 1} \geqslant 1 - 2/63 ( p + 1
  )$. Thus there exists a constant $d > 0$ such that $\alpha_{p - 1} \geqslant 1 -
  d/p$.                                                     
  
  {\tmstrong{III.}} Let $\eta = \min ( c, d )$, then $\alpha_{\ast}
  \leqslant 1 - \eta/p$. By Theorem 7 we have

  \[ \| K_t ( x, . ) - \pi \|_{\tmop{TV}} \leqslant \frac{1}{2}
  \pi_{\ast}^{- \frac{1}{2}} \alpha_{\ast}^t \leqslant \frac{1}{2}
  \pi_{\ast}^{- \frac{1}{2}} ( 1 - \frac{\eta}{p} )^t \]

  where $\pi_{\ast} = \min_x \{ \pi ( x ) \} = \pi ( c_0 ) = 1/p^2$.
  Thus,
  \[ \| K_t ( x, . ) - \pi \|_{\tmop{TV}} \leqslant \frac{p ( 1
  - \frac{\eta}{p} )^t}{2}. \]
  Therefore, $ \| K_t ( x, . ) - \pi \|_{\tmop{TV}} \leqslant 1/2 e
  $ if $p ( 1 - \eta/p )^t/2 \leqslant 1/2 e$. It is equivalent to $1
  + p + t \ln ( 1 - \eta/p ) \leqslant 0$ or $1 + p + t ( -
  \eta/p - \eta^2/p^2 - \ldots ) \leqslant 0$. So we can choose any $t$ such that $\eta ( 1
  + \ln p ) p \leqslant t$. This concludes the proof.        
\end{proof}

\subsection{Probabilistic approach}

In this section, we will prove a better bound for $\tau_p$ by using
probabilistic approach, the coupling technique.

\begin{proposition}
  Given $p \in \mathcal{P}, p \equiv 3$ (mod $4$). If $( 1 -
  p^2 ( p - 1 )/( 1 + p )^4 )^n < 1/2 e $ then $r_p
  \leqslant 4 n$.
\end{proposition}

\begin{proof}
  We create two copies of random walk $c_1^n$. The first one starts from $c_i$
  for fixed $i$ and the second one starts from circles with distribution $\pi$. In
  step $m$, suppose that we are in circle $c_t$ in the first walk and in circle $c_s$ in      the second walk for some $s, t$. If $t = s$ then in the next step, we choose the step
  of the second walk which is the same with the first's. Otherwise, let they
  walk by $A_1$ independently. It is clearly that both random walks have
  the same Markov kernel $A_1$ and the second one has the distribution $\pi$.
  
  {\tmem{Claim}}: $A_1^4 ( c_i, c_j ) \geqslant \frac{p^2 ( p - 1 )}{( 1 + p
  )^4} \pi ( c_j )$ for all $i, j$.
  
  Proof. We consider 4 separated cases
  
  1) $i = j = 0$. There are $|C_1 |^4 = ( p + 1 )^4$ possible ways to go by 4 steps. We
  first go by any 2 steps. In the last two, we just go backward then it is clear that we go
  back to the starting point. Therefore, at least $|C_1 |^2 = ( p + 1 )^2$
  ways to go from $c_0$ to $c_0 $. It implies that
 \[
  A_1^4 ( c_0, c_0 ) \geqslant \frac{(
  p + 1 )^2}{( p + 1 )^4} > \frac{p^{^2} ( p - 1 )}{( 1 + p )^4} \pi ( c_j ) =
  \frac{p - 1}{( p + 1 )^4}.
 \]  
  2) $i = 0, j \neq 0$. We have
  \[
  A_1^4 ( c_0, c_j ) = \sum_{l, k} A_1 ( c_0, c_1 ) A_1 ( c_1, c_l )
  A_1 ( c_l, c_k ) A_1 ( c_k, c_j ).
  \]
  By Theorem 2, we have for each $l \neq 0$ then exists $k \neq 0$ such that 
  \[A_1 ( c_l, c_k ) A_1 ( c_k, c_j )> 0. \]
  
  But  $A_1 ( c_u, c_v ) > 0$, then $A_1 ( c_u, c_v ) \geqslant
  1/(p + 1)$. Since $\Pr ( l = 0 ) = n_{11}^0 = 1/(p + 1)$ we have  
  \[ A_1^4 ( c_0, c_j ) \geqslant \frac{1}{( p + 1 )^2}
  \Pr ( l \neq 0 ) = \frac{p}{( p + 1 )^3}.\]  
  
  Therefore, we have
  \[A_1^4 ( c_0, c_j ) \geqslant \frac{p^3}{( p + 1 )^4}
  \frac{( p + 1 )}{p^2} > \frac{p^2 ( p - 1 )}{( p + 1 )^4} \pi ( c_j ).\]
  
  3) $i \neq 0, j = 0$. Similar as 2), we have 
  \[A_1^4 ( c_i, c_0 ) \geqslant \frac{p^3}{( p + 1 )^4}
  \frac{( p + 1 )}{p^2} > \frac{p^3}{( p + 1 )^4} \frac{1}{p^2} > \frac{p^2 (
  p - 1 )}{( p + 1 )^4} \pi ( c_j ).\]
  
  4) $i, j \neq 0$. We have 
  \[ A_1^4 ( c_i, c_j ) = \sum_{t, l, k} A_1 ( c_i, c_t ) A_1 ( c_t, c_l
  ) A_1 ( c_l, c_k ) A_1 ( c_k, c_j ). \]
  Similar as in 2), we have 
  \[A_1^4 (c_i, c_j ) \geqslant \frac{1}{( p + 1 )^2} \Pr ( l \neq 0 ) = \frac{1}{( p +
  1 )^2} ( 1 - \Pr ( l = 0 ) ) .\]  
  But \[\Pr ( l = 0 ) = \sum_t A_1 ( c_i, c_t ) A_1 ( c_t, c_0 ) = n_{11}^i
  \leqslant \frac{2}{p + 1}.\]
  So we have \[A_1^4 ( c_i, c_j ) \geqslant \frac{p - 1}{( p
  + 1 )^3} = \frac{p^2 ( p - 1 )}{( p + 1 )^4} \pi ( c_j ).\]  
  This finishes the proof of the claim.                                     
  
  Set $c = 1 - p^2 ( p - 1 )/( 1 + p )^4$. From the claim, we have
  \begin{align*}  
  d_{\tmop{TV}} ( A_1^4 ( c_i, . ), \pi ) &= \frac{1}{2} \sum_j |
  \pi ( c_j ) - A_1^4 ( c_i, c_j ) | \\
  &= \sum_{j : A_1^4 ( c_i, c_j ) < \pi ( c_j
  )} ( \pi ( c_j ) - A_1^4 ( c_i, c_j ) )\\  
  &\leqslant \sum_{j : A_1^4 ( c_i, c_j ) < \pi
  ( c_j )} \pi ( c_j ) ( 1 - c ) \\
  &\leqslant 1 - c.
  \end{align*}
                          
  Applying the above inequality $t$ times, we have
  
  \[d_{\tmop{TV}} ( A_1^{4 t} ( c_i, . ), \pi ) \leqslant ( 1 - c
  )^t.\]
  Therefore, if $( 1 - p^2 ( p - 1 )/( 1 + p )^4 )^n < 1/2
  e$ then  $d_{\tmop{TV}} ( A_1^{4 t} ( c_i, . ), \pi ) \leqslant 1/2
  e$ and $\tau_p \leqslant 4 n$ (Q.E.D)
\end{proof}

As a corollary we have the following theorem.

\begin{theorem}
  $\tau_p = O ( p )$.
\end{theorem}

\begin{proof}
  The inequality $( 1 - p^2 ( p - 1 )/( 1 + p )^4)^n <
  1/2 e $ is equivalent to
  \[n \log \left( 1 - \frac{p^2 ( p - 1 )}{( 1 + p)^4} \right) < - \ln 2 - 1. \]
Therefore $n$ need to satisfy
  
\[n \left( \frac{p^2 ( p - 1 )}{( p + 1 )^4} + \frac{p^4 (
  p - 1 )^2}{( 1 - p )^8} + \ldots \right) > 1 + \ln 2.\]
Thus, $n > \frac{( 1 + \log 2 ) ( p + 1 )^4}{p^2 ( p - 1 )}$. This conludes the proof.
\end{proof}

\begin{center}
\Large\bfseries{Acknowlegment}
\end{center}
The original impetus of this research was given by my supervisor, A.Prof Norman Wildberger.
I am grateful for his helpful discussions and numerous invaluable suggestions.

\end{document}